# Reliability of first order numerical schemes for solving shallow water system over abrupt topography

T. Morales de Luna[a,*], M.J. Castro Díaz[b], C. Parés[b]

[a]*Dpto. de Matemáticas. Universidad de Córdoba. Campus de Rabanales. 14071 Córdoba. Spain.*
[b]*Dpto. de Análisis Matemático. Facultad de Ciencias. Universidad de Málga. Campus de Teatinos, s/n. 29071 Málaga, Spain.*

**Abstract**

We compare some first order well-balanced numerical schemes for shallow water system with special interest in applications where there are abrupt variations of the topography. We show that the space step required to obtain a prescribed error depends on the method. Moreover, the solutions given by the numerical scheme can be significantly different if not enough space resolution is used. We shall pay special attention to the well-known hydrostatic reconstruction technique where it is shown that the effect of large bottom discontinuities might be missed and a modification is proposed to avoid this problem.

## 1. Introduction

The numerical analysis of the shallow water system has been extensively studied in recent years. This system may be written in the form:

[*]Corresponding author
  *Email addresses:* tomas.morales@uco.es (T. Morales de Luna),
castro@anamat.cie.uma.es (M.J. Castro Díaz), pares@anamat.cie.uma.es (C. Parés)

$$\begin{cases} \dfrac{\partial h}{\partial t} + \dfrac{\partial q}{\partial x} = 0, \\ \dfrac{\partial q}{\partial t} + \dfrac{\partial}{\partial x}\left(\dfrac{q^2}{h} + \dfrac{g}{2}h^2\right) = gh\dfrac{\partial H}{\partial x}. \end{cases} \quad (1.1)$$

The variable $x$ makes reference to the axis of the channel and $t$ is time; $q(x,t)$ and $h(x,t)$ represent the mass-flow and the thickness, respectively; $g$, the acceleration due to gravity; $H(x)$, the depth measured from a fixed level of reference; $q(x,t) = h(x,t)u(x,t)$, with $u$ the depth averaged horizontal velocity.

This model has the structure of a system of balance laws:

$$\partial_t w + \partial_x F(w) = S(w)\partial_x H, \quad (1.2)$$

with

$$w = \begin{pmatrix} h \\ q \end{pmatrix}, \quad F(w) = \begin{pmatrix} q \\ \dfrac{q^2}{h} + \dfrac{g}{2}h^2 \end{pmatrix}, \quad S(w) = \begin{pmatrix} 0 \\ gh \end{pmatrix}.$$

By considering $H$ as an (artificial) unknown whose value is given by the initial conditions, the system can also be written as a first order quasilinear system:

$$W_t + \mathcal{A}(W) \cdot W_x = 0, \quad x \in \mathbb{R}, t > 0, \quad (1.3)$$

with $W = (h, q, H)^t$ and

$$\mathcal{A}(W) = \begin{pmatrix} 0 & 1 & 0 \\ -u^2 + c^2 & 2u & -c^2 \\ 0 & 0 & 0 \end{pmatrix}, \quad (1.4)$$

where

$$u = \dfrac{q}{h}, \quad c = \sqrt{gh}.$$

The eigenvalues of $\mathcal{A}(W)$ are

$$\lambda_1 = u - c, \quad \lambda_2 = u + c, \quad \lambda_3 = 0,$$

and the system is strictly hyperbolic provided that $h > 0$ and $|u| \neq c$.

The Riemann invariants corresponding to the null eigenvalues are:

$$q = constant, \quad h + \dfrac{q^2}{2gh^2} - H = constant. \quad (1.5)$$

In the particular case $q = 0$ these curves are straight lines in $h, q, H$ space:

$$q = 0, \qquad h - H = constant. \tag{1.6}$$

Once the depth function $H(x)$ has been fixed, the stationary solutions of the system are given by (1.5). In particular, those given by (1.6) correspond to water at rest solutions.

Models of the form (1.2) or (1.3) also appear in many other physical problems. It is well known that standard methods that solve correctly systems of conservation laws can fail in the presence of source terms, especially when approaching equilibria or near to equilibria solutions. In the context of shallow water equations Bermúdez and Vázquez-Cendón introduced in [27] the condition called conservation property or C-property: a scheme is said to satisfy this condition if it solves correctly the steady state solutions corresponding to water at rest. This idea of constructing numerical schemes that preserve some equilibria, which are called in general well-balanced schemes, has been studied by many authors (see, for instance, [16], [6], [21] [11], [19], [2], [7], [3], [26], [15], [12] among many others).

This study was motivated by the following observation: given a fixed mesh, different well-balanced numerical methods based on different treatments of the source terms can lead to very different numerical solutions. Figure 1 illustrates this phenomenon: a numerical test consisting of a shallow layer of fluid flowing up a ramp has been considered, and the water elevation obtained with different first order numerical schemes using the same mesh are compared. Some similar test cases as well as the considered numerical schemes will be described in detail below. As it can be seen, very different behaviors of the fluid layer is predicted by the considered numerical schemes: in some cases the fluid goes up the ramp in a smooth way, in some others a strong shock going to the left is predicted. It is worth noting that these huge differences are not related to a convergence problem: these differences can be done arbitrarily small by refining the mesh and/or by increasing the order of accuracy of the method. Nevertheless, for different first order methods, the refinement level required to obtain a solution close to the common limit can be significantly bigger than those necessary for some others. Moreover, when solving these problems with a first order finite volume method, one approaches the solution with a piece-wise constant function. This means that in practice we have a *discontinuity* of the bottom at each interface, which motivates the study not only for steep topography but in the presence of a discontinuous bottom, where we shall see that this kind of phenomenon arises. This phenomenon implies a difficulty: in practical applications one usually has a given mesh whose refinement level depends on the available bathymetric data and on the available computational power. Therefore, it is important to know how reliable the numerical solutions are if there are abrupt variations of the topography and the mesh is not fine enough. Roughly speaking, a numerical treatment of the source term is expected to be more reliable than others if the refinement level of the mesh which is required to obtain a numerical solution whose behavior is close enough to the limit is lower.

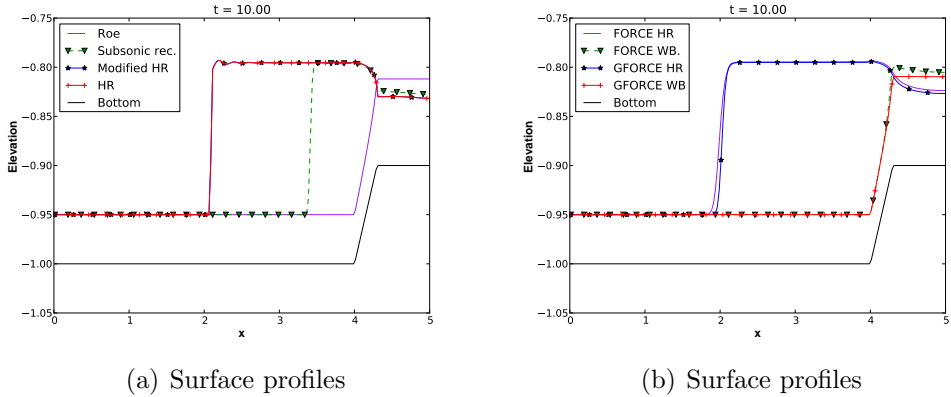

(a) Surface profiles  (b) Surface profiles

Figure 1: Comparison of numerical solutions obtained with different schemes

The goal of this article is to run a battery of numerical tests in order to compare different numerical schemes based on different treatment of the source terms and to draw conclusions about their reliability. The different numerical treatment of the source terms considered here belong to two families: on the one hand, we consider numerical methods based on an upwinding treatment of the source term. This treatment has been introduced by Roe for scalar problems [22] and by Bermúdez and Vázquez-Cendón in [27]. On the other hand, we consider numerical methods based on the hydrostatic reconstruction technique introduced in [2] and some extensions and modifications. A phenomenon which is closely related to the one that motivated this paper was observed in [14] for numerical methods based on the hydrostatic reconstruction technique: it was seen that the numerical solutions obtained on a fixed mesh with a first order numerical scheme for a test consisting of a shallow layer of fluid going down a ramp were independent of the slope beyond a critical threshold, what is not in good agreement with the physics of the problem. This difficulty was also discussed in [5]. A modification of the hydrostatic reconstruction technique will be introduced here that partially overcomes this difficulty.

Although the numerical analysis of system (1.1) with discontinuous depth function $H$ is not addressed here, the theoretical and numerical difficulties arising in this case can give some insight about the origin of the phenomenon to be studied here. When the depth function $H$ is discontinuous, the source term cannot be defined within the distributional framework: it is a non-conservative product that may be defined in many different ways: see, for instance, [13]. Depending on the chosen definition of the source term across a discontinuity of $H$, one has many different ways to define what is a weak solution of the problem. According to the more natural and extended definition, the simple waves associated to a bottom discontinuity are contact discontinuity preserving the Riemann invariants of the third characteristic field (1.5): see [16], [1], [17] for example. Nevertheless, some alternative definitions of weak solution can be found in the literature: see [4], [23] for example.

Once the notion of weak solution has been fitted, a concept of entropy solution is necessary to discard the unphysical ones. In the case of the shallow water system, there is a known entropy pair that allows one to define such a concept: a weak solution is said to be an entropy solution if the inequality

$$\partial_t \widetilde{\eta}(w, H) + \partial_x \widetilde{G}(w, H) \leq 0 \qquad (1.7)$$

is satisfied in the distribution sense, where

$$\begin{aligned} \eta(w) &= h\frac{u^2}{2} + \frac{g}{2}h^2, & G(w) &= \left(\frac{u^2}{2} + gh\right)hu, \\ \widetilde{\eta}(w, H) &= \eta(w) - ghH, & \widetilde{G}(w, H) &= G(w) - ghuH. \end{aligned} \qquad (1.8)$$

Nevertheless, in this case the addition of an entropy notion does not ensure the uniqueness of solution: in the case of resonance problems, i.e. when one of the eigenvalues of the matrix system changes its sign: see, for instance, [17].

Let us remark that, although a discontinuous depth function $H$ is out of the theoretical assumptions used to derive shallow water models (one of them is that the normal vector to the bottom is almost parallel to the vertical axis), it is relevant to consider here the possibility of a discontinuous $H$ for the following reasons:

- The family of numerical methods described below are first order numerical schemes based on a piecewise constant approximation of the bottom function. Due to this, the solution of a Riemann problem with discontinuous values of $H$ has to be approached at every intercell. Even for a continuous bottom function, a step is found at discrete level and this step can be large for steep bottom if the mesh is not fine enough. As different approximate Riemann solvers may converge to different solutions, these local approximations may be very different to each other: this is the reason of the different behaviors observed in Figure 1.

  From the numerical point of view, the main difficulty of solving the system with a discontinuous depth function comes from the fact that the solutions provided by two different numerical methods can converge to functions which are weak solutions according to different definitions, or to different entropy solutions according to the same definition of weak solutions.

- Even though bottom discontinuities are out of the theoretical range of validity of the model, the shallow water model still provides reasonable results in the presence of small enough bottom steps if the weak solutions are allowed to have stationary contact discontinuities at the bottom jumps.

This paper is organized as follows: in Section 2 and for the sake of completeness, the different numerical schemes to be compared are briefly presented. Section 3 is devoted to the hydrostatic reconstruction, focusing on the difficulty mentioned before and a modified version of the hydrostatic reconstruction technique is proposed to overcome it. Finally, Section 4 will present different numerical simulations comparing previously cited numerical schemes.

## 2. Upwind treatment of the source term

We consider as usual a set of computing cells $I_i = [x_{i-1/2}, x_{i+1/2}], i \in \mathbb{Z}$. We shall assume, for the sake of simplicity, that these cells have a constant size $\triangle x$ and that $x_{i+1/2} = i\triangle x$. $x_i = (i - 1/2)\triangle x$ is the center of the cell $I_i$. Let $\triangle t$ be the time step and $t^n = n\triangle t$.

The following notation will be used for the approximation of the cell averages of the exact solution:

$$w_i^n = \begin{bmatrix} h_i^n \\ q_i^n \end{bmatrix}, \quad W_i^n = \begin{bmatrix} w_i^n \\ H_i \end{bmatrix}. \tag{2.1}$$

We consider numerical schemes that can be written in the form

$$w_i^{n+1} = w_i^n - \frac{\Delta t}{\Delta x}(F_{i+1/2} - F_{i-1/2}) + \frac{\Delta t}{\Delta x}(S_{i-1/2}^+ + S_{i+1/2}^-), \tag{2.2}$$

where

$$F_{i+1/2} = \mathcal{F}(w_i^n, w_{i+1}^n), \tag{2.3}$$

$\mathcal{F}$ being a continuous function such that

$$\mathcal{F}(w,w) = F(w), \tag{2.4}$$

and

$$S_{i+1/2}^\pm = \mathcal{S}^\pm(W_i^n, W_{i+1}^n), \tag{2.5}$$

where $\mathcal{S}^-$ and $\mathcal{S}^+$ have two components: a centered approximation of the source term and the upwind treatment.

Let us introduce some examples:

2.1. Roe method

The extension of the Roe method to (1.2) is given by (2.2) with the choices

$$F_{i+1/2} = \frac{1}{2}\left(F(w_i^n) + F(w_{i+1}^n)\right) - \frac{1}{2}|J_{i+1/2}|(w_{i+1}^n - w_i^n), \tag{2.6}$$

$$S_{i+1/2}^\pm = \mathcal{P}_{i+1/2}^\pm S_{i+1/2}(H_{i+1} - H_i), \tag{2.7}$$

where

$$J_{i+1/2} = \begin{bmatrix} 0 & 1 \\ -(u_{i+1/2}^n)^2 + (c_{i+1/2}^n)^2 & 2u_{i+1/2}^n \end{bmatrix} \tag{2.8}$$

with

$$u_{i+1/2}^n = \frac{\sqrt{h_i^n}u_i^n + \sqrt{h_{i+1}^n}u_{i+1}^n}{\sqrt{h_i^n} + \sqrt{h_{i+1}^n}}, \quad c_{i+1/2}^n = \sqrt{g\frac{h_i^n + h_{i+1}^n}{2}}. \tag{2.9}$$

$S_{i+1/2}$ is given by:

$$S_{i+1/2} = \begin{bmatrix} 0 \\ (c_{i+1/2}^n)^2 \end{bmatrix}, \tag{2.10}$$

and
$$\mathcal{P}^{\pm}_{i+1/2} = \frac{1}{2}\left(Id \pm |J_{i+1/2}|J^{-1}_{i+1/2}\right). \tag{2.11}$$

Here $Id$ represents the identity matrix and

$$|J_{i+1/2}| = K_{i+1/2}|\Lambda_{i+1/2}|K^{-1}_{i+1/2}, \tag{2.12}$$

where $|\Lambda_{i+1/2}|$ is the diagonal matrix whose coefficients are the absolute value of the eigenvalues of $J_{i+1/2}$ and $K_{i+1/2}$ is the matrix whose columns are associated eigenvectors. This numerical scheme is exactly well-balanced for water-at-rest solutions and is second order accurate for general stationary solutions: see [27], [21].

2.2. FORCE and GFORCE methods

FORCE and GFORCE schemes were introduced in [24] and [25] for conservative systems. Their well-balanced extensions to nonconservative systems were introduced in [10]. For the particular case of the shallow water system, these extensions can be written in the form (1.2) with

$$\begin{aligned} F_{i+1/2} = &\frac{1}{2}\left(F(w_i^n) + F(w_{i+1}^n)\right) \\ &- \frac{1}{2}\left((1-\omega)\frac{\Delta x}{\Delta t}Id + \omega\frac{\Delta t}{\Delta x}J^2_{i+1/2},\right)(w_{i+1}^n - w_i^n), \end{aligned} \tag{2.13}$$

$$\begin{aligned} S^{\pm}_{i+1/2} = &S_{i+1/2}(H_{i+1} - H_i) \\ &\pm \left((1-\omega)\frac{\Delta x}{\Delta t}J^{-1}_{i+1/2} + \omega\frac{\Delta t}{\Delta x}J_{i+1/2},\right)S_{i+1/2}(H_{i+1} - H_i) \end{aligned} \tag{2.14}$$

where $J_{i+1/2}$ and $S_{i+1/2}$ are given by (2.8) and (2.10) respectively. FORCE and GFORCE methods correspond, respectively, to $\omega = 1/2$ and $\omega = 1/(1+CFL)$, where $CFL \in (0,1]$ is the CFL number. Notice that the choices $\omega = 0$ and $\omega = 1$ correspond to extensions of Lax-Friedrichs and Lax-Wendroff methods. This numerical scheme is also exactly well-balanced for water at rest solution and second order accurate for general stationary solutions.

Remark that (2.14) is not well defined when one of the eigenvalues of $J_{i+1/2}$ vanishes. This means that when there is a change from subcritical to supercritical regime, this scheme is not properly defined. To correct this, we may replace $J_{i+1/2}^{-1}$ by

$$\widetilde{J}_{i+1/2}^{-1} = \frac{1}{\mu} \begin{pmatrix} 0 & 1 \\ (c_{i+1/2}^n)^2 & 2u_{i+1/2} \end{pmatrix}, \quad (2.15)$$

where $\mu = \max\{\varepsilon, |1 - Fr^2|\} \cdot \text{sgn}(1 - Fr^2)$, $Fr^2 = u_{i+1/2}^2/(c_{i+1/2})^2$ is the square of Froude number and $\varepsilon > 0$ is some given small value.

Another possibility would be to use the following modification of the discrete source term:

$$S_{i+1/2}^{\pm} = S_{i+1/2}(H_{i+1} - H_i)$$
$$\pm \left( (1-\omega)\frac{\Delta x}{\Delta t}(J_{i+1/2}^*)^{-1} + \omega \frac{\Delta t}{\Delta x} J_{i+1/2}, \right) S_{i+1/2}(H_{i+1} - H_i) \quad (2.16)$$

where:

$$J_{i+1/2}^* = \begin{bmatrix} 0 & 1 \\ (c_{i+1/2}^n)^2 & 0 \end{bmatrix}. \quad (2.17)$$

This modified method is still exactly well-balanced for water-at-rest solutions but only first order accurate for general stationary solutions. We refer to [10] for further details.

2.3. *Path-conservativity*

The upwind numerical methods introduced in this section satisfy the following equalities:

$$S_{i+1/2}^{\pm} = 0 \text{ if } W_i^n = W_{i+1}^n, \quad (2.18)$$

and

$$S_{i+1/2}^+ + S_{i+1/2}^- = \int_0^1 S(\Psi_w(s; W_i^n, W_{i+1}^n))\frac{\partial \Psi_H}{\partial s}(s; W_i^n, W_{i+1}^n)ds, \quad (2.19)$$

where

$$\Psi(s; W_i^n, W_{i+1}^n) = W_i^n + s\left(W_{i+1}^n - W_i^n\right), \quad (2.20)$$

with $\Psi_w$ and $\Psi_H$ the two first and last component of $\Psi$ respectively. These properties imply that these methods can be interpreted as path-conservative numerical methods in the sense defined in [20] for (1.3), that is, these methods are formally consistent with a particular definition of weak solution of the Riemann problem with a discontinuous depth function $H$: more precisely, they are consistent with the definition of nonconservative product provided by the theory developed in [13] corresponding to the choice of the family of straight segments (2.20). The choice of the family of straight segments is a natural choice to obtain well-balanced numerical methods for water-at-rest solutions, as the equation (1.6) corresponds to straight lines in the $h, q, H$ space: [18] for details.

## 3. Hydrostatic reconstruction technique

*3.1. The original formulation*

The hydrostatic reconstruction technique was introduced in [2]. It presents a straightforward method for obtaining a well-balanced scheme for the shallow water system based on a numerical flux for the homogeneous system ($H = constant$). Moreover, this technique has the advantage of preserving the *good properties* of the homogeneous solver such as positivity of water thickness and entropy inequalities. The idea is to consider a consistent numerical flux $\mathcal{F}(w_l, w_r)$ for the homogeneous shallow water system, that is (1.1) with $H$ constant, and the corresponding numerical scheme

$$w_i^{n+1} = w_i^n + \frac{\Delta t}{\Delta x}(F_{i-1/2} - F_{i+1/2}), \quad F_{i+1/2} = \mathcal{F}(w_i, w_{i+1}). \quad (3.1)$$

Then define at each intercell the values

$$H_{i+1/2} = \min(H_i, H_{i+1/2}), \quad q_{i+1/2}^- = h_{i+1/2}^- u_i, \quad q_{i+1/2}^+ = h_{i+1/2}^+ u_{i+1}, \quad (3.2)$$

$$\begin{aligned} h_{i+1/2}^+ &= (h_{i+1} - H_{i+1} + H_{i+1/2})_+, \\ h_{i+1/2}^- &= (h_i - H_i + H_{i+1/2})_+, \end{aligned} \quad (3.3)$$

and define the scheme for non-homogeneous shallow water equations in the form (2.2) with

$$F_{i+1/2} = \mathcal{F}(w_{i+1/2}^-, w_{i+1/2}^+) \quad (3.4)$$

$$S_{i+1/2}^- = \begin{pmatrix} 0 \\ p(h_i) - p(h_{i+1/2}^-) \end{pmatrix}, \quad S_{i+1/2}^+ \begin{pmatrix} 0 \\ p(h_{i+1}) - p(h_{i+1/2}^+) \end{pmatrix}, \quad (3.5)$$

where $p(h) = \frac{g}{2}h^2$.

The resulting numerical method is non-negative and semi-discrete entropy satisfying provided that the homogeneous flux also satisfies these properties.

*3.2. Path-conservativity*

The resulting scheme can be also interpreted as a path-conservative method for (1.3), i.e. (2.18)-(2.19) are again satisfied but for a family of paths different from the straight segments (2.20). The family of paths in now defined as follows: given $W_i^n$ and $W_{i+1}^n$, the path linking them is composed by:

1. the segment on the space $h, q, H$ defined by the equations:

$$q = hu_i^n, \quad h - H = h_i^n - H_i,$$

   linking the states $W_i^n$ and

$$W_{i+1/2}^- = \begin{bmatrix} h_{i+1/2}^- \\ u_i^n h_{i+1/2}^- \\ H_{i+1/2} \end{bmatrix};$$

2. the segment linking the states $W_{i+1/2}^-$ and

$$W_{i+1/2}^+ = \begin{bmatrix} h_{i+1/2}^+ \\ u_{i+1}^n h_{i+1/2}^+ \\ H_{i+1/2} \end{bmatrix};$$

3. the segment on the space $h, q, H$ defined by the equations:

$$q = hu_{i+1}^n, \quad h - H = h_{i+1}^n - H_{i+1},$$

   linking the states $W_{i+1/2}^+$ and $W_{i+1}^n$.

Notice that when $H_{i+1/2} = H_i$ the first segment degenerates to a point and when $H_{i+1/2} = H_{i+1}$ is the third one that degenerates to a point. Figure 3(a) shows the third segment in the case $H_{i+1/2} = H_i$.

The discretized source terms (3.5) are the integral of the source term

$$\int_C gh \, dH, \tag{3.6}$$

computed by using the first and third segments as the integration path $C$. Observe that the second one does not contribute to the source term as it lies in a plane $H = constant$.

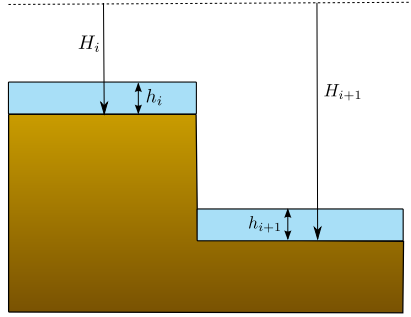

Figure 2: Case $h^+_{i+1/2} = 0$

When $H_{i+1/2} = H_i$ and $h^n_{i+1} - H_{i+1} + H_i < 0$ or $H_{i+1/2} = H_{i+1}$ and $h^n_i - H_i + H_{i+1} < 0$, the paths are modified to grant the positivity of the numerical scheme. Let us consider the first case, the second one being analogous: the situation is that shown in Figure 2. In this case, $h^+_{i+1/2} = 0$, and the path linking the states is composed by:

1. The segment linking $W^n_i$ and

$$W^+_{i+1/2} = \begin{bmatrix} 0 \\ 0 \\ H_i \end{bmatrix};$$

2. the arc of the graph

$$h = (h_{i+1} - H_{i+1} - H)_+, \quad q = hu^n_{i+1}$$

linking the states $W^+_{i+1/2}$ and $W^n_{i+1}$.

Figure 3(b) shows the second piece of this path (dashed line).

Observe that, in this case, the only piece of the path that contributes to the discretized source term is the segment linking the state $W^n_{i+1}$ and $[0, 0, H_{i+1} - h^n_{i+1}]^T$ in which the straight line

$$h - H = h^n_{i+1} - H_{i+1}, \quad q = hu^n_{i+1}$$

intersects the axis

$$h = 0, \quad q = hu^n_{i+1}.$$

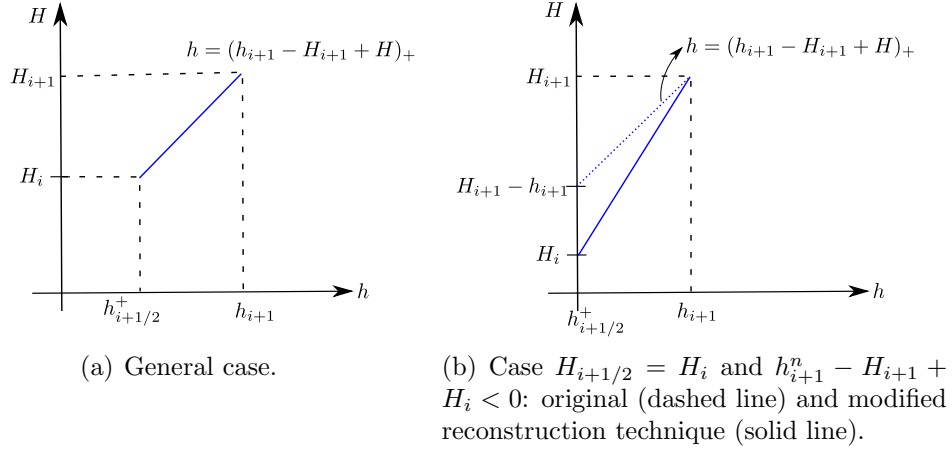

(a) General case.　　(b) Case $H_{i+1/2} = H_i$ and $h_{i+1}^n - H_{i+1} + H_i < 0$: original (dashed line) and modified reconstruction technique (solid line).

Figure 3: Hydrostatic reconstruction technique.

Therefore, the numerical solution does not depend on the value of $H_i$ in this case: this is the reason of the phenomenon observed in [14] and mentioned in the Introduction.

The interpretation of the hydrostatic reconstruction technique in terms of the choice of a particular family of paths has been on the basis of the generalization of the technique introduced in [11] that allows to obtain a first order numerical method which is exactly well-balanced for every stationary solution. It is also on the basis of the simpler modification proposed in the next subsection to partially overcome this loss of accuracy in the presence of big steps.

*3.3. Modified hydrostatic reconstruction*

In order to take into account the amplitude of the step, we propose a modification of the hydrostatic reconstruction technique. In terms of choice of the family of paths, the proposal is applied whenever $H_{i+1/2} = H_i$ and $h_{i+1}^n - H_{i+1} + H_i < 0$; or $H_{i+1/2} = H_{i+1}$ and $h_i^n - H_i + H_{i+1} < 0$. In the first case, the idea is to replace the path described above by:

1. The segment linking $W_i^n$ and

$$W_{i+1/2}^+ = \begin{bmatrix} 0 \\ 0 \\ H_i \end{bmatrix};$$

2. the segment linking $W^+_{i+1/2}$ and $W^n_{i+1}$.

(See the solid line in Figure 3(b)). The analogous modification is proposed when $H_{i+1/2} = H_{i+1}$ and $h^n_i - H_i + H_{i+1} < 0$.

In order to keep a similar notation to the one used in the original hydrostatic reconstruction, the numerical scheme is written in the form (2.2) where now:

$$F_{i+1/2} = \mathcal{F}(w^-_{i+1/2}, w^+_{i+1/2}) \tag{3.7}$$

$$S^-_{i+1/2} = \begin{pmatrix} 0 \\ p(h_i) - p(h^-_{i+1/2}) + T^-_{i+1/2} \end{pmatrix},$$
$$S^+_{i+1/2} = \begin{pmatrix} 0 \\ p(h_{i+1}) - p(h^+_{i+1/2}) + T^+_{i+1/2} \end{pmatrix}. \tag{3.8}$$

Some easy computations of the integral (3.6) through the modified paths lead to the following results:

1. General case.

   If $h_i - H_i + \min(H_i, H_{i+1}) \geq 0$ and $h_{i+1} - H_{i+1} + \min(H_i, H_{i+1}) \geq 0$, then the standard hydrostatic reconstruction is used. Therefore:

   $$T^\pm_{i+1/2} = 0.$$

2. Large steps at the intercell.

   If $h_i - H_i + \min(H_i, H_{i+1}) < 0$ or $h_{i+1} - H_{i+1} + \min(H_i, H_{i+1}) < 0$, then the path is modified

   $$T^-_{i+1/2} = p(h^-_{i+1/2}) - p(h_i) + g\frac{h_i + h^-_{i+1/2}}{2}(H_{i+1/2} - H_i), \tag{3.9}$$

   $$T^+_{i+1/2} = p(h^+_{i+1/2}) - p(h_{i+1}) + g\frac{h_{i+1} + h^+_{i+1/2}}{2}(H_{i+1} - H_{i+1/2}). \tag{3.10}$$

The resulting numerical method takes into account the total height of the step and it is well-balanced for water at rest solutions, provided that no wet/dry fronts are present, as in this case, the method coincides with the original one. Nevertheless, this modification can produce unphysical solutions in the presence of an emerging bottom: see Figure 4.

Consider for example the case of emerging bottom to the right (Figure 4(a)) and assume that water is at rest $u = 0$. It can be easily checked that this steady state is not preserved by the modified hydrostatic reconstruction. On the other hand, if $u_i > 0$ in this situation, two different cases will be distinguished:

- The mechanical energy at the wet cell is enough to make the fluid to go up the step. Then, the full step should be considered when computing the numerical fluxes.

- There is not enough energy at the wet cell, and thus the step acts as an obstacle.

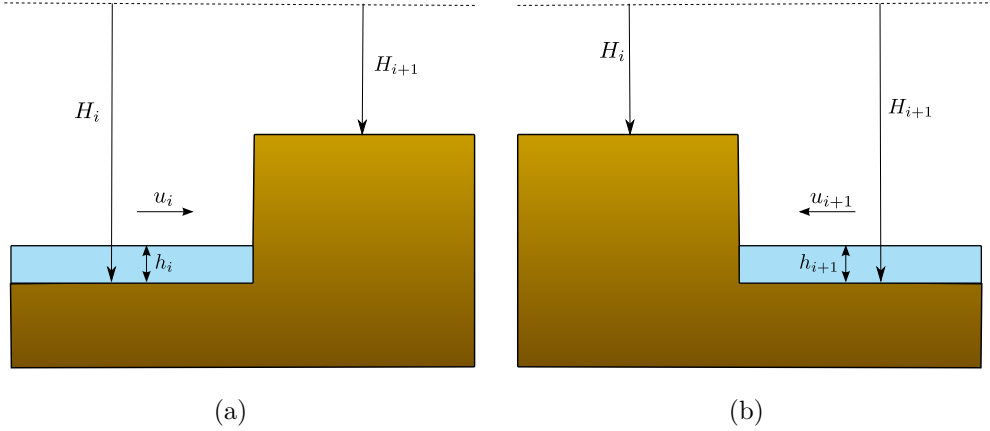

Figure 4: Emerging bottom

Following the ideas presented in [8], the fluid should go up the step in the case of emerging bottom to the right and $u_i > 0$ if

$$\frac{u_i^2}{2} + g(h_i - H_i + H_{i+1}) > \frac{3}{2}\sqrt{(gh_i u_i)^3}, \qquad (3.11)$$

and in the case of emerging bottom to the left and $u_{i+1} < 0$ if

$$\frac{u_{i+1}^2}{2} + g(h_{i+1} - H_{i+1} + H_i) > \frac{3}{2}\sqrt{(gh_{i+1} u_{i+1})^3} \qquad (3.12)$$

Therefore, in the presence of emerging bottom situations, the numerical method is applied as follows:

- In water at rest situations or if the mechanical energy is not enough to make the water go up the step, the original hydrostatic reconstruction technique is applied.

- If the mechanical energy is enough, the modified reconstruction technique is applied to take into account the total height of the step.

The non-negativity-preserving property of the numerical methods based on the hydrostatic reconstruction technique is also satisfied for the modified technique:

**Proposition 3.1.** *Consider a numerical scheme in the form (2.2), (3.7),(3.8) with the modified hydrostatic reconstruction described above. If $\mathcal{F}$ preserves the non-negativity of $h$ by interface, then so does the modified hydrostatic reconstruction.*

The proof of this result is essentially the same as the one presented for hydrostatic reconstruction in [6].
The entropy preserving property is discussed in Appendix A.

*3.4. Subsonic reconstruction scheme*

The subsonic reconstruction scheme was introduced in [7] as a generalization of the hydrostatic reconstruction. This solver intends to preserve a further class of steady state solutions for the shallow water system while preserving the good properties of being positive and entropy satisfying. The definition of this scheme is rather technical and for the sake of brevity it is not included here. We refer the reader to [7] for a complete description.

## 4. Numerical simulations

A battery of numerical tests have been designed to compare the numerical solutions obtained with the different discretizations of the source terms introduced in Sections 2 and 3. The original or the modified reconstruction technique are used by default combined with Roe's flux for the homogeneous problem, but in some cases, the FORCE or the GFORCE numerical fluxes are considered. The following labels are used in the Figures:

- Roe: the numerical method given by (2.2), (2.6), (2.7).

- HR: the hydrostatic reconstruction technique combined with Roe's numerical flux for the homogeneous problem.

- Modified HR: the modified hydrostatic reconstruction technique combined with Roe's numerical flux for the homogeneous problem.

- FORCE HR: the hydrostatic reconstruction technique combined with FORCE numerical flux for the homogeneous problem.

- FORCE WB: the numerical method given by (2.2), (2.13), (2.14), and $\omega = 1/2$.

- GFORCE HR: the hydrostatic reconstruction technique combined with GFORCE numerical flux for the homogeneous problem.

- GFORCE WB: the numerical method given by (2.2), (2.13), (2.14), and $\omega = 1/(1 + CFL)$.

- Subsonic reconstruction: the numerical method described in [7].

### 4.1. Test 1

This first test intends to show the problem observed in [14]. To do so, let us reproduce the test shown there by considering in the interval $[0,3]$ the initial condition

$$h(x, t=0) = 0.02, \quad q(x, t=0) = 0.01, \quad H(x) = \frac{\alpha}{100}x, \qquad (4.1)$$

for different values of slope $\alpha$. Boundary conditions to the left are given by $q(x=0, t) = 0.01$, $h(x=0, t) = 0.02$, and free boundary conditions are set on the right boundary.

As it was remarked by Delestre, water height differs up to a given critical value of the slope. For instance, in Figure 5(a) we see that when we use 50 points in the interval $[0,3]$, HR scheme reproduces always the same water thickness. This is related to what has already been discussed before: for the resolution of 50 points, we found at each intercell at discrete level a jump which is large enough for HR scheme to neglect it. The modified HR scheme corrects this fault as it is shown in Figure 5(b). This can be corrected also by using a second order reconstruction technique as it was shown in [14].

We remark that here the bottom is continuous, so one expects that when $\Delta x \to 0$ this problem will not be found as we see in Figure 6.

### 4.2. Test 2

This test shows that all the numerical methods give similar result when the depth function is continuous and the steps of the piecewise constant discretization of the bottom are small enough. We consider a test taken

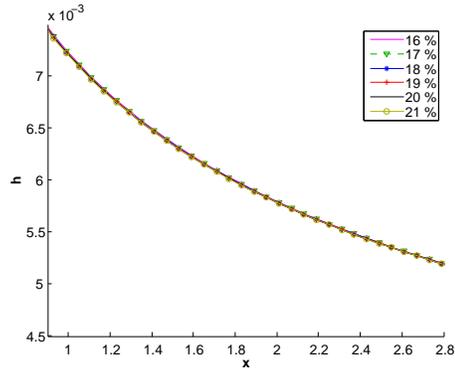
(a) HR scheme

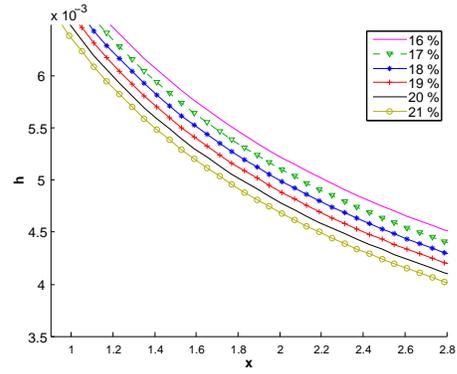
(b) Modified HR

Figure 5: Test 1: water height for 50 points for $\alpha = 16, 17, \ldots, 21$

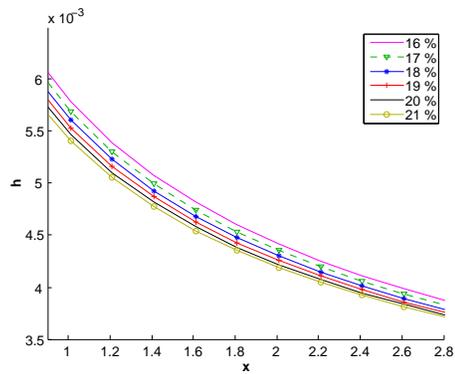
(a) HR scheme

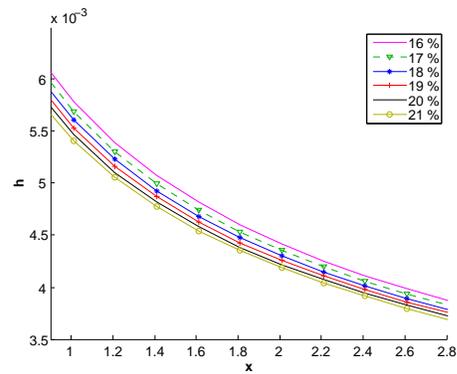
(b) Modified HR

Figure 6: Test 1: water height for 150 points for $\alpha = 16, 17, \ldots, 21$

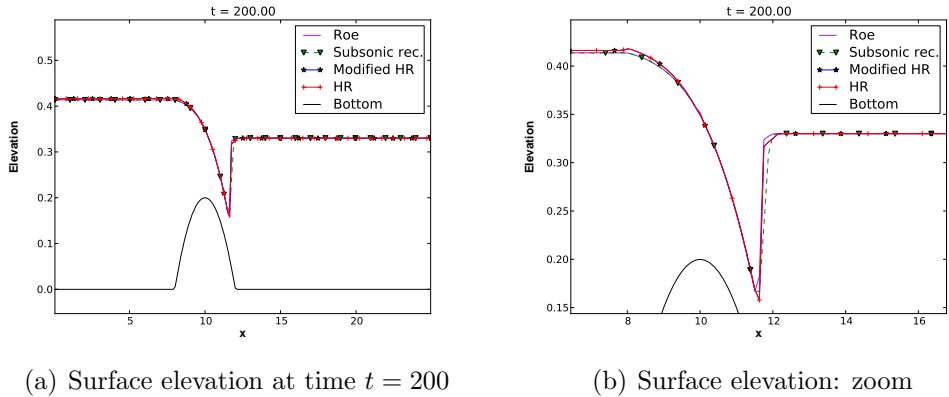

(a) Surface elevation at time $t = 200$

(b) Surface elevation: zoom

Figure 7: Test 2: transcritical solution

from [6] consisting of a transcritical flow with a shock over a bump. The interval is $[0, 25]$, the initial condition and depth bottom is given by

$$h(x, t = 0) = 0.33, \quad q(x, t = 0) = 0.18,$$
$$H(x) = \begin{cases} -0.2 + 0.05(x - 10)^2, & \text{if } 8 < x < 12, \\ 0, & \text{otherwise}, \end{cases} \quad (4.2)$$

and the boundary conditions are given by $q(x = 0, t) = 0.18$, $h(x = 25, t) = 0.33$.

We consider a uniform mesh composed by 200 cells. The height of the steps of the discretized bottom are lower than $1/200$. The numerical results are shown in Figure 7. As expected, all the numerical results are similar.

4.3. Test 3

The goal of this test is to show that, in the presence of a discontinuity of the bottom function, the numerical solutions provided by different numerical methods may converge to different limits which are weak solutions of the problem according to different notions, or to different entropy solutions according to a same notion. Let us remark again that, even if discontinuous bottoms are out of the range of validity of the model, the interest of such a test comes from the fact that numerical methods based on piecewise discontinuous approximation of the depth functions, as it is the case for the ones considered here, have to deal with bottom steps at every interface. And the local differences of the approximate Riemann solutions may lead to very different global behavior for not fine enough meshes.

This phenomenon is studied in the following test: we consider the following initial condition and depth function

$$h(x, t=0) = 0.1, \quad q(x, t=0) = 0.15, \quad H(x) = \begin{cases} 0.1, & \text{for } x < 0.5, \\ H_r, & \text{otherwise,} \end{cases} \quad (4.3)$$

with increasing values of $H_r$ from 0.1 to 0.5. The left boundary conditions are given by $h(x=0,t) = 0.1$, $q(x=0,t) = 0.1$. and open boundary conditions are considered at $x = 1$. A uniform mesh of 200 cells is considered. Figure 8 shows the water elevation obtained for $H_r = 0.15$ and $H_r = 0.45$ at different times: the differences of the numerical solutions provided by the considered numerical methods are apparent.

In order to systematically study these differences, we compute the values of the water thickness $h_r$ to the right of the discontinuity obtained with the different numerical methods when a steady state is reached. They are also compared with the value $h_r$ such that the Riemann invariants corresponding to the states $W_l = (0.1, 0.1, 0.1)$ and $W_r = (h_r, 0.1, H_r)$ coincide: as it was mentioned in the Introduction, according to the most used definition of weak solution, the exact stationary solution of the problem will be given by a contact discontinuity linking the states $W_l$ and $W_r$. Therefore, this value of $h_r$ is labelled as 'exact' in the Figures.

The results are shown in Figure 9. As it has been reported previously, the results given by the hydrostatic reconstruction technique are the same once the height of the step are greater than a critical value. The modification introduced here corrects this problem. Nevertheless, we remark that none of the schemes reproduce the exact solution, specially in presence of big discontinuities of the bottom. This is a known fact that has been discussed in [9].

### 4.4. Test 4

We consider again a step-like bottom, but in this case the water layer flows from the bottom to the top of the step. The initial condition is given again by (4.3) but the depth function is now:

$$H(x) = \begin{cases} 0.8, & \text{for } x < 0.5 \\ H_r, & \text{otherwise ,} \end{cases} \quad (4.4)$$

with $H_r$ varying from 0.8 to 0.4. A uniform mesh of 200 cells is again considered. The boundary conditions $h(x=0,t) = 0.1$, $q(x=0,t) = 0.15$ are imposed at $x = 0$ and open boundary conditions are considered at $x = 1$. Different numerical solutions are shown in Figures 10 to 12: as it can be seen, the bigger the step, the bigger the differences. In order to compare the different behaviors, $F_r^2 = u^2/(gh)$ is also depicted.

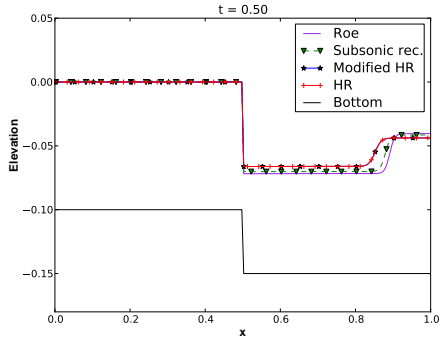
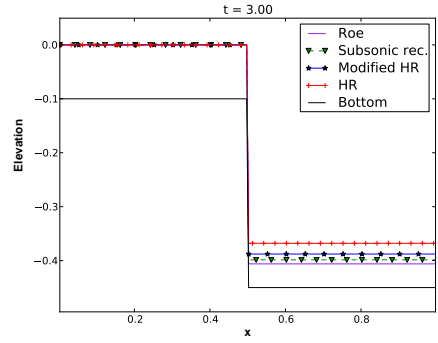

(a) $H_r = 0.15$, $t = 0.50$

(b) $H_r = 0.45$, $t = 3$

Figure 8: Test 3: solutions obtained for two different values of $H_r$ in (4.3)

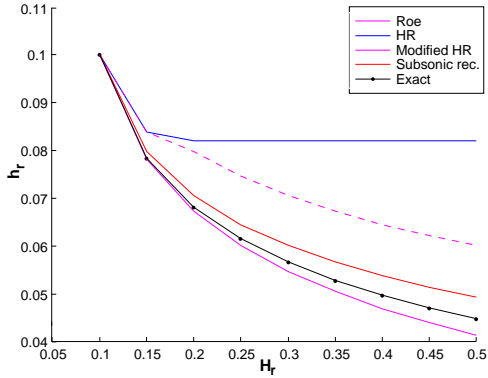
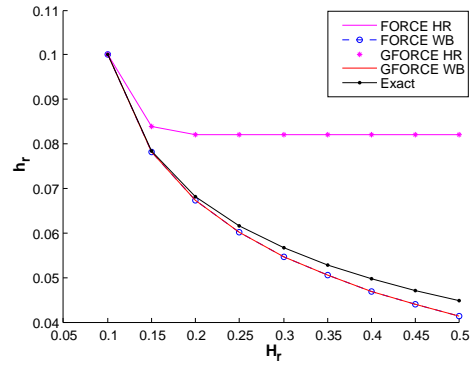

(a)

(b)

Figure 9: Test 3: water thickness $h_r$ to the right of the step versus $H_r$, when a steady-state is reached

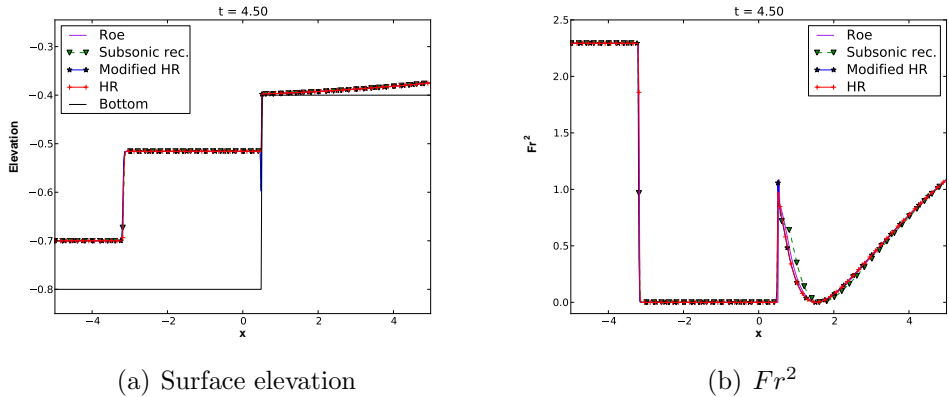

(a) Surface elevation  (b) $Fr^2$

Figure 10: Test 4: Surface elevation and Froude number (squared) versus $x$ for given value $H_r = 0.4$, $t = 4.5$

As in the previous test, we compute the values of the thickness $h_l$ and $h_r$ to the left and to the right of the step once a steady state has been reached. The results are shown in Figure 13 and Figure 14. Both hydrostatic reconstructions give essentially the same results but they are very different from the ones obtained with the subsonic reconstruction scheme or Roe. Remark that the value $h_r$ is the same for all the schemes for $H_r = 0.8$, which is the case corresponding to a flat bottom. The same occurs for $H_r \leq 0.5$ when all the numerical solutions predict a reflected shock.

In this case, we have transcritical regions and (2.15) should be used. We remark that results shown in Figure 14(b) depend on the parameter $\varepsilon$. Moreover, in this particular case, if we use (2.16), the results provided by FORCE WB and GFORCE WB are closer to those given by their HR version. Moreover, given that the bottom is discontinuous and we obtain transcritical solutions, the differences could be related to the non-uniqueness of solution for this case. This will be discussed in more detail in following example.

*4.5. Test 5*

We consider now a test consisting of a shallow layer of water flowing up a ramp with increasing slope. More explicitly, we consider the following depth

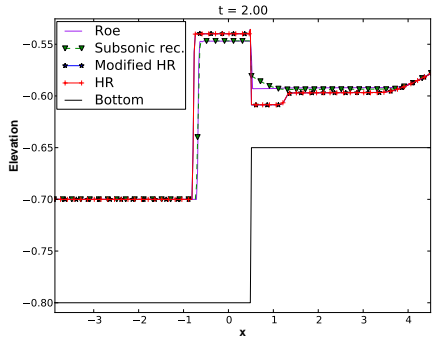 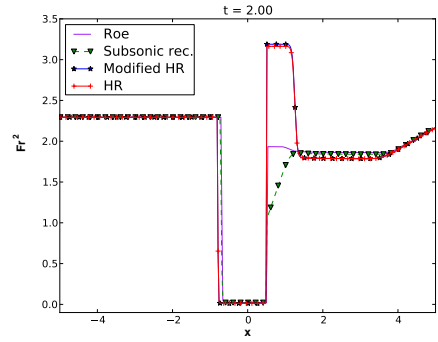

(a) Surface elevation  (b) $Fr^2$

Figure 11: Test 4: Surface elevation and Froude number (squared) versus $x$ for given value $H_r = 0.65$, $t = 2$

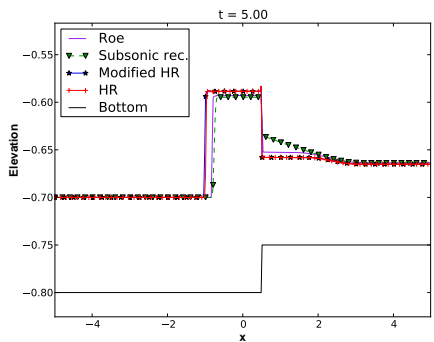 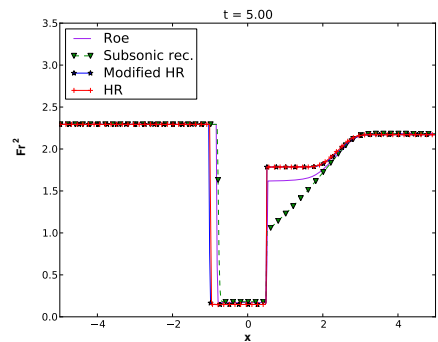

(a) Surface elevation  (b) $Fr^2$

Figure 12: Test 4: Surface elevation and Froude number (squared) versus $x$ for given value $H_r = 0.75$, $t = 5$)

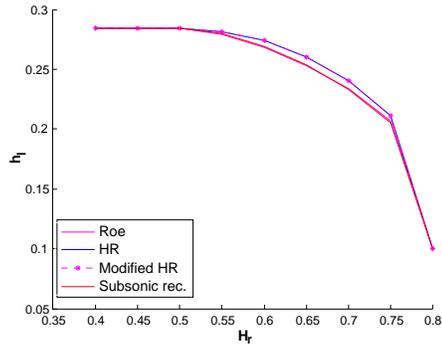
(a) Left state $h_l$

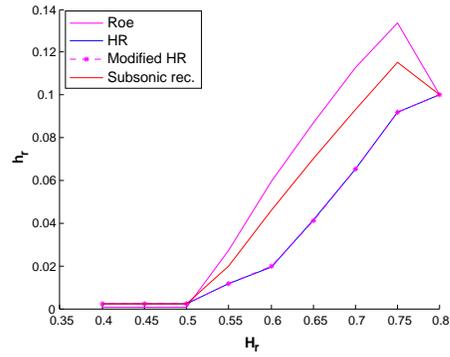
(b) Right state $h_r$

Figure 13: Test 4: Comparison of the water thickness $h_l$ (respectively $h_r$) to the left (respectively to the right) of the step versus $H_r$ once a steady state has been reached

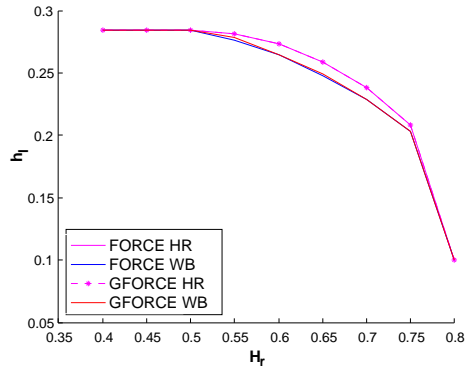
(a) Left state $h_l$

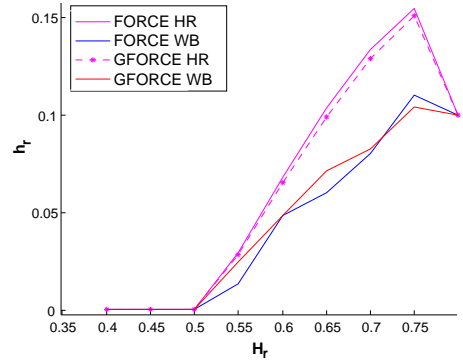
(b) Right state $h_r$

Figure 14: Test 4: Comparison of the water thickness $h_l$ (respectively $h_r$) to the left (respectively to the right) of the step versus $H_r$ once a steady state has been reached

function and initial condition:

$$H(x) = \begin{cases} 1, & \text{if } x < x_l, \\ 1 + \dfrac{2.8}{4 - x_l}(x - x_l), & \text{if } x_l \leq x < 4, \\ 0.2, & \text{if } x > 4 \end{cases} \quad (4.5)$$

$$h(x, t = 0) = (H - 0.9)_+, \quad q(x, t = 0) = \begin{cases} 0.9, & \text{if } h(x, t = 0) > 0, \\ 0, & \text{elsewhere} \end{cases}, \quad (4.6)$$

with boundary conditions $h(x = 0, t) = 0.1$, $q(x = 0, t) = 0.9$, and $x_l < 4$ being a given parameter. Figure 15 shows the bottom and the initial surface elevation for $x_l = 3.75$.

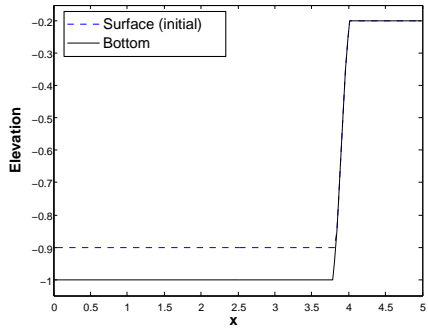

Figure 15: Test 5: example of initial condition

Figures 16 to 21 show the computed water thickness $h_l$ (respectively $h_r$) before (respectively after) the rise in the topography at time $t = 2.5$ for three uniform meshes of the interval $[0, 5]$ composed by 200, 400, and 800 cells.

Remark that, as the slope becomes steeper, the values obtained with both hydrostatic reconstructions are further from those obtained for the other schemes. The bigger number of cells, the lesser the differences, as one would expect. Moreover, if the grid is not fine enough, extremely different behaviors can be observed to the left of the slope: either the water goes up the ramp in a smooth way or a shock traveling to the left develops: see Figure 22. This is related to what has already been mentioned in previous tests: for steep

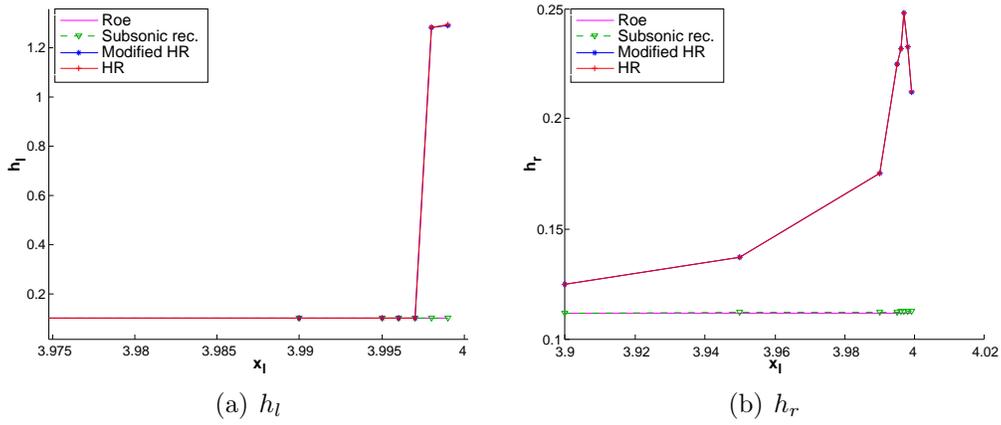

Figure 16: Test 5: Computed water thickness $h_l$ (respectively $h_r$) before (respectively after) the rise in the topography versus the initial point of the slope $x_l$ using 200 points

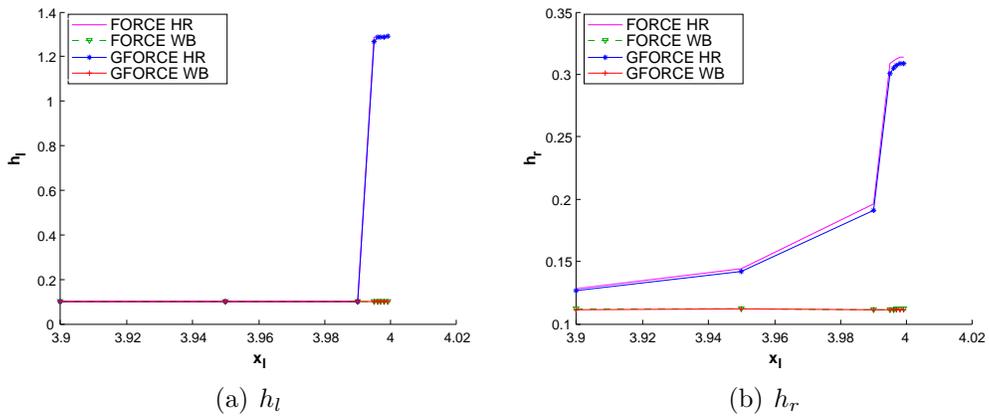

Figure 17: Test 5: Computed water thickness $h_l$ (respectively $h_r$) before (respectively after) the rise in the topography versus the initial point of the slope $x_l$ using 200 points

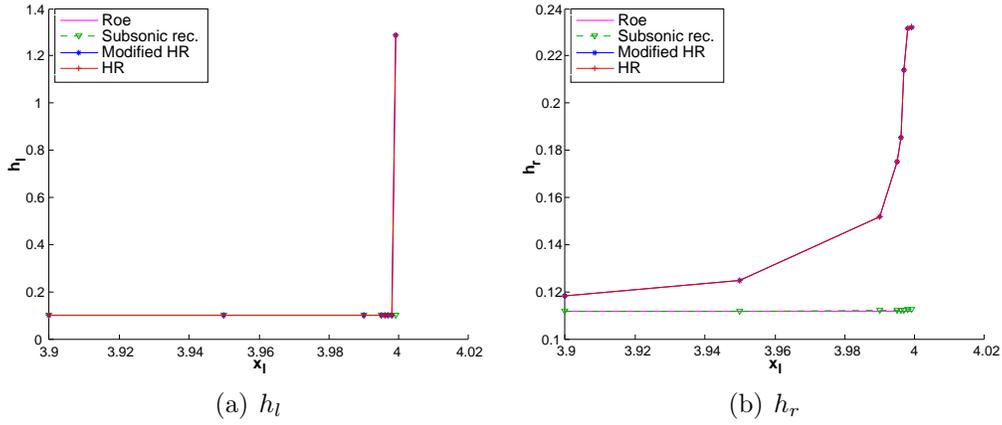

Figure 18: Test 5: Computed water thickness $h_l$ (respectively $h_r$) before (respectively after) the rise in the topography versus the initial point of the slope $x_l$ using 400 points

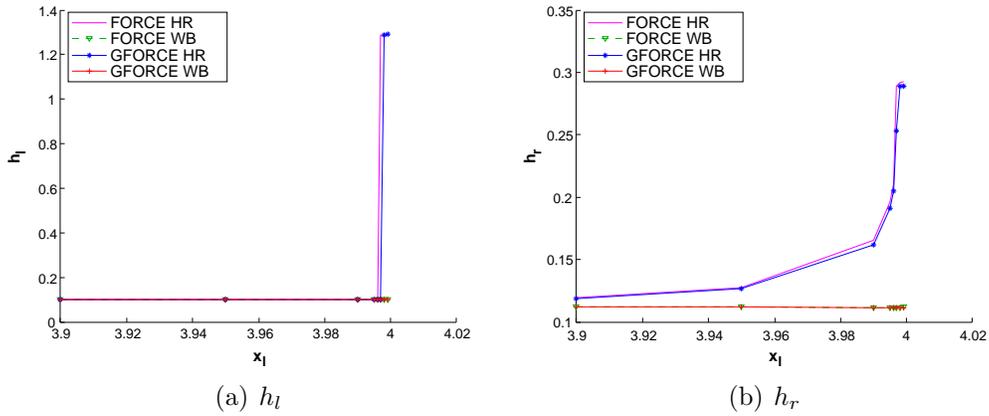

Figure 19: Test 5: Computed water thickness $h_l$ (respectively $h_r$) before (respectively after) the rise in the topography versus the initial point of the slope $x_l$ using 400 points

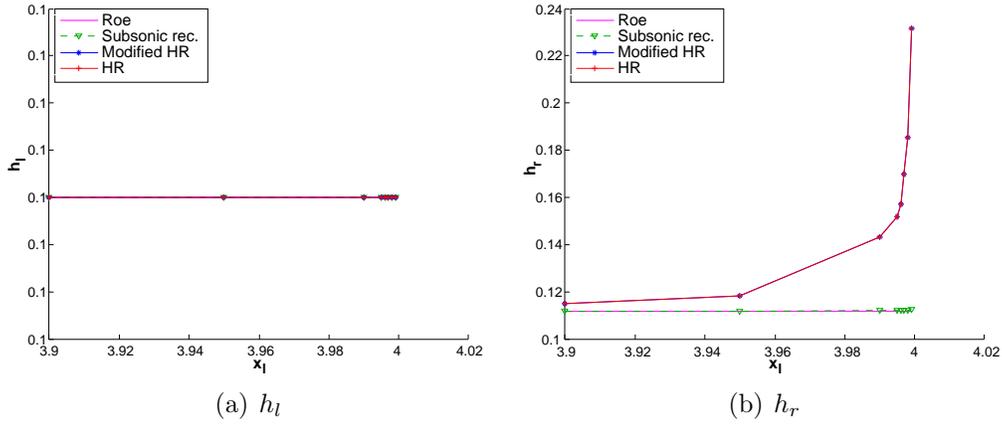

Figure 20: Test 5: Computed water thickness $h_l$ (respectively $h_r$) before (respectively after) the rise in the topography versus the initial point of the slope $x_l$ using 800 points

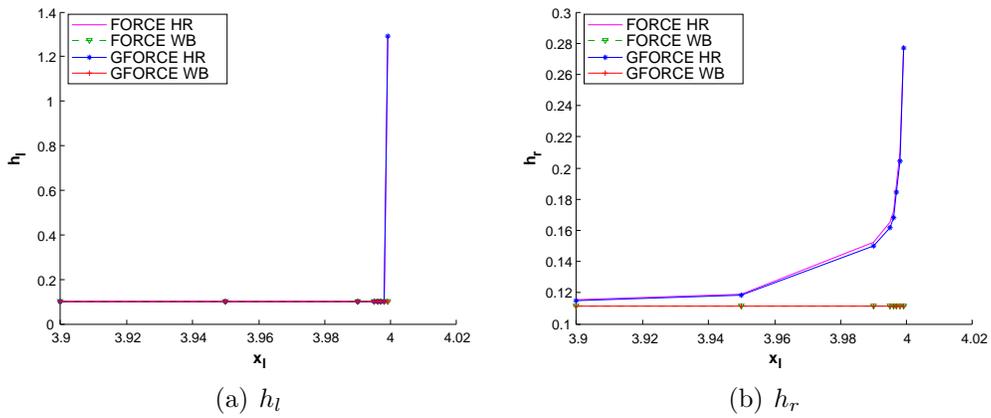

Figure 21: Test 5: Computed water thickness $h_l$ (respectively $h_r$) before (respectively after) the rise in the topography versus the initial point of the slope $x_l$ using 800 points

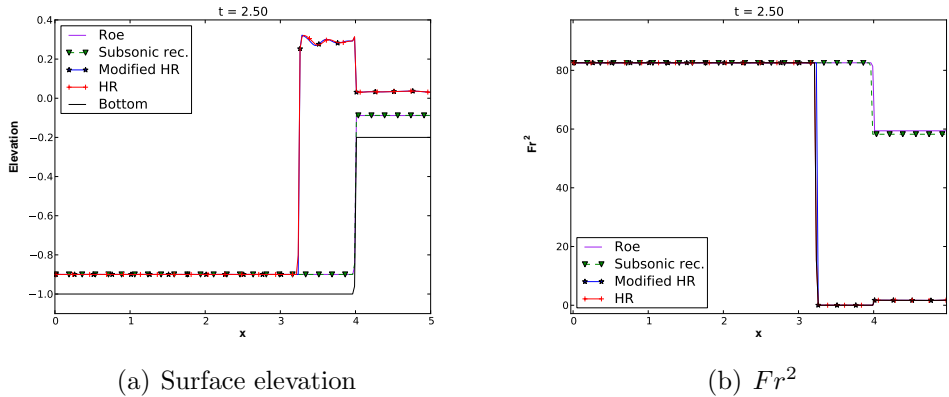

(a) Surface elevation  (b) $Fr^2$

Figure 22: Test 5: Surface elevation and Froude number (squared) versus $x$ for a given value $x_l$ close to 4 at $t = 2.5$

bottom and not enough number of cells the numerical scheme "sees" a "large" discontinuity for the bottom. Moreover, in this case the bottom topography is continuous and all the numerical solutions converge to the same supercritical stationary solution. But if the mesh is not fine enough, the slope is treated as a bottom step and the numerical solutions provided by the different numerical methods may be very different: see, for instance, Figure 22. In Figure 22(b) the Froude numbers corresponding to the numerical solutions are drawn and it may be seen that, while some of the numerical solutions are supercritical, some others involve transcritical shocks. This fact is very likely related to the non-uniqueness of weak solution in resonant situations.

In general, we see that the schemes behave in a similar way for small slopes and large number of cells. But for steep slopes, on the one hand, we need to use more points for those based on hydrostatic reconstruction; on the other hand, Roe, subsonic reconstruction, FORCE WB and GFORCE WB behave in a similar way even for small number of cells. We remark again that we have here a continuous bottom and the numerical solutions should converge to a supercritical steady-state solution. So, the differences among them are not related to a problem of convergence, but rather a problem of resolution or reliability of the schemes if we do not use enough number of points.

*4.6. Test 6*

We consider finally a test similar to the previous one, consisting of a shallow layer of water flowing down a ramp with increasing slope. More explicitly, we consider the following depth function and initial condition

$$h(x, t = 0) = 0.5 \quad q(x, t = 0) = 1.2, \tag{4.7}$$

$$H(x) = \begin{cases} 0.1, \text{ if } x \leq 0.2, \\ 0.1 + \dfrac{H_r - 0.1}{x_r - 0.2}(x - 0.2), \text{ if } 0.2 < x \leq x_r, \\ H_r \text{ if } x > x_r \end{cases} \tag{4.8}$$

with boundary conditions $h(x = 0, t) = 0.5$, $q(x = 0, t) = 1.2$, $x_r = 0.2 + \Delta l$, $H_r = 0.1 + \Delta H$ with $\Delta l$ and $\Delta H$ given parameters.

The smooth stationary solution reached can be computed from (1.5) and it is shown, together with the initial condition, in Figure 23 for a particular value of $\Delta l$ and $\Delta H$.

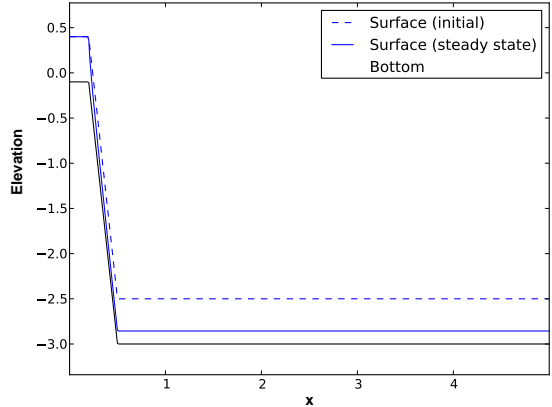

Figure 23: Test 6: Surface elevations corresponding to the initial condition and to the stationary solution.

The numerical methods are applied to this test with different values of $\Delta l$ and $\Delta H$ and different number of cells in the interval $[0, 5]$. More explicitly, 8 uniform meshes having 100, 200, 400, 800, 1600, 3200, 6400, and 12800 cells have been considered. Then, the exact stationary solution is computed for every value of $\Delta l$ and $\Delta H$ and it is compared with the ones computed with the numerical schemes. As expected, all the numerical stationary solutions converge to the exact one as the mesh is refined. Nevertheless, the number of cells required to satisfy a given error bound strongly depends on the numerical method: Figures 24 to 26 show the number of cells needed in order to obtain a $L^1$-error lower than 0.008. Incomplete graphs mean that the prescribed error bound is not reached for the finer grid of 12800 cells.

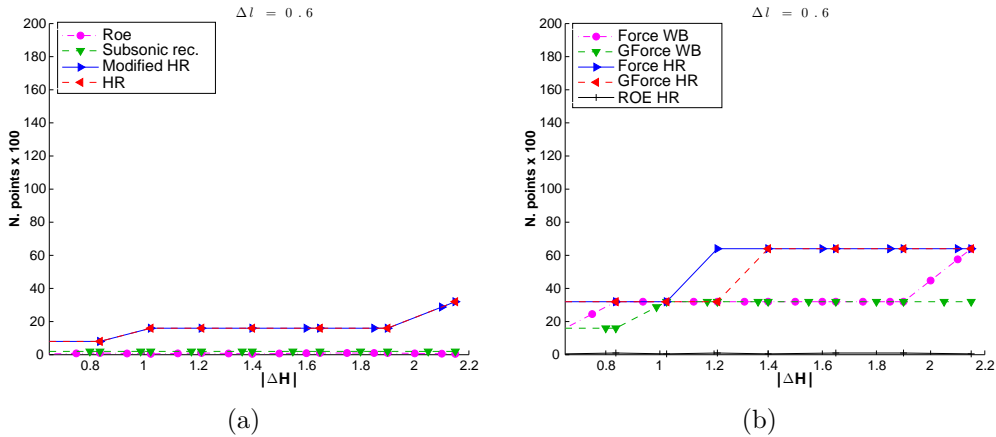

Figure 24: Test 6: Number of cells needed to obtain $L^1$-error $\leq 0.008$ versus size of the step $\Delta H$

As it can be seen, Roe or the subsonic reconstruction schemes need only a small number of cells to obtain the prescribed bound, even for the steepest ramp. But this is not the case for the other methods: this number dramatically increases as the slope of the ramp increases.

## 5. Conclusions

We have compared different well-balanced numerical schemes for the shallow water system based on different discretizations of the source term in test cases in which either the bottom is discontinuous or has steep slopes. Even though a discontinuous depth function is out of the theoretical assumptions used to derive shallow water models, the interest of such study is motivated by the fact that numerical schemes shown here are first order schemes based on piecewise constant approximations. Due to this, the solution of a Riemann problem with discontinuous values of the depth function has to be approached at each intercell. While the numerical solutions can converge to different limits for discontinuous bottom, this is not the case for continuous depth functions. Nevertheless, the space step required to obtain an error which is lower than a prescribed bound dramatically depends on the method. Moreover, if the mesh is not fine enough, the numerical solutions can exhibit completely different behaviors. A numerical method is considered to be more reliable as the space step required to achieve a prescribed accuracy is bigger.

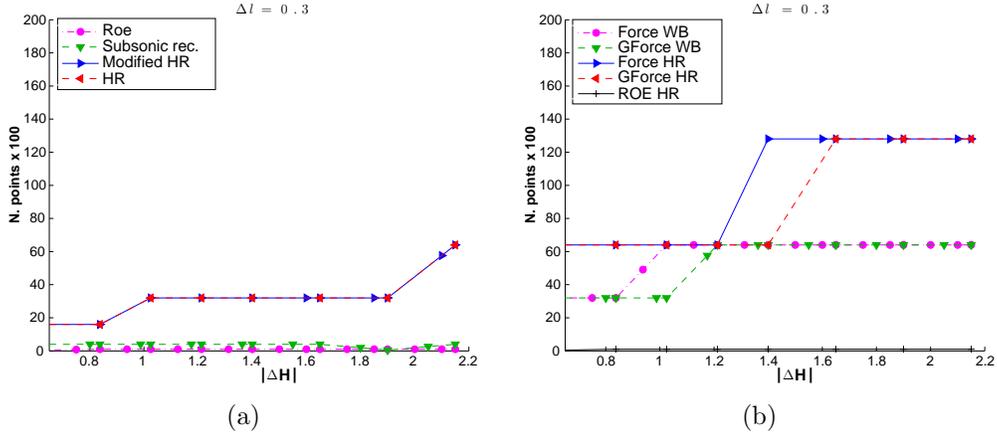

Figure 25: Test 6: Number of cells needed to obtain $L^1$-error $\leq 0.008$ versus size of the step $\Delta H$

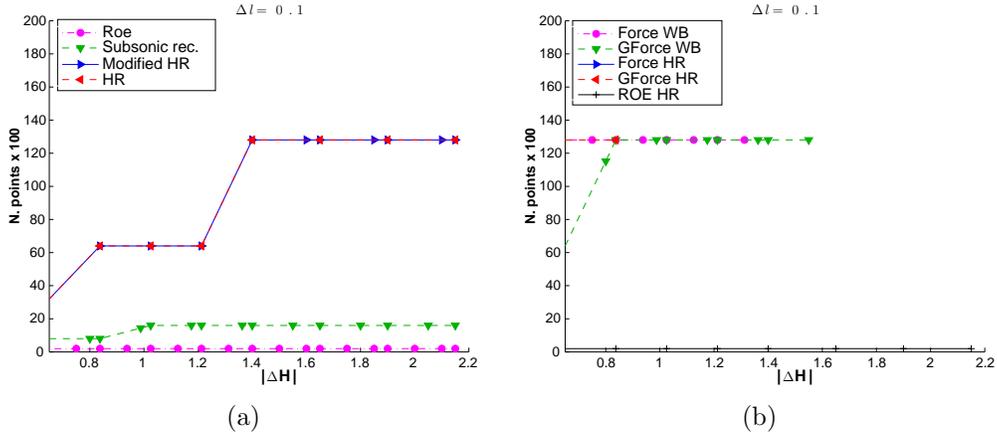

Figure 26: Test 6: Number of cells needed to obtain $L^1$-error $\leq 0.008$ versus size of the step $\Delta H$

We have considered two different approaches to discretize the source term: an upwind treatment or the hydrostatic reconstruction technique. While this latter technique, introduced in [2], allows one to easily define a robust scheme for shallow water equations, its reliability (in the sense mentioned above) is not very high: this is due to the fact that the effect of large bottom discontinuities might be missed. We have introduced a modification of this technique that takes into account the total height of the bottom discontinuities. Nevertheless, the numerical tests show that both the original and the modified hydrostatic reconstruction are less reliable than other numerical methods such as the subsonic reconstruction scheme [7] or the Roe scheme [21]. In general, it has been observed that, for continuous steep bottoms, a large number of discretization points in the domain is needed for numerical methods based on the hydrostatic reconstruction techniques, FORCE or GFORCE methods, while Roe or the subsonic reconstruction schemes give good accuracy even with a small number of cells.

## A. Appendix

In [2] it was shown that whenever $\mathcal{F}$ is an entropy satisfying numerical flux for the homogeneous shallow water system, the hydrostatic reconstruction scheme satisfies a semi-discrete entropy inequality for the full system in the sense given by the following definition:

**Definition A.1.** *A numerical scheme for (1.2) in the form (2.2) - (2.3) is said to satisfy a semi-discrete entropy inequality with respect to to the entropy-entropy flux pair $(\widetilde{\eta}, \widetilde{G})$ given by (1.8) if there exists a numerical entropy flux function $\mathcal{G}(W_l, W_r)$ consistent with the exact entropy flux such that for any $W_l, W_r$ we have*

$$\widetilde{G}(W_r) + \nabla_w \widetilde{\eta}(W_r) \left( \mathcal{F}_r(W_l, W_r) - F(w_r) \right) \leq \mathcal{G}(W_l, W_r)$$
$$\leq \widetilde{G}(W_l) + \nabla_w \widetilde{\eta}(W_l) \left( \mathcal{F}_l(W_l, W_r) - F(w_l) \right). \quad (A.1)$$

In most situations, the modified hydrostatic reconstruction scheme presented here coincides with the original one. Thus, in most cases we will obtain an entropy satisfying scheme.

To study the general case, we shall use a result presented in [7]:

**Proposition A.2.** *Consider a numerical scheme in the form*

$$w_i^{n+1} = w_i^n - \frac{\Delta t}{\Delta x}(F_{i+1/2}^- - F_{i-1/2}^+), \quad (A.2)$$

*with $F_{i+1/2}^\pm = \mathcal{F}^\pm(W_i, W_{i+1})$ and*

$$\mathcal{F}^-(W_l, W_r) = \mathcal{F}(w_l^*, w_r^*) + \begin{pmatrix} 0 \\ p(h_l) - p(h_l^*) + T^-(W_l, W_r) \end{pmatrix},$$
$$\mathcal{F}^+(W_l, W_r) = \mathcal{F}(w_l^*, w_r^*) + \begin{pmatrix} 0 \\ p(h_r) - p(h_r^*) + T^+(W_l, W_r) \end{pmatrix}, \quad (A.3)$$

*where $\mathcal{F}$ stands for a given consistent numerical flux for the homogeneous shallow water system that verifies a semi-discrete entropy inequality for the entropy pair $(\eta, G)$ in (1.8) and the interface values $w_l^*, w_r^*$ are derived from a local reconstruction procedure. Denote by $\mathcal{F}(w_l^*, w_r^*) = (\mathcal{F}^h, \mathcal{F}^q)$.*

*A sufficient condition for the scheme to be semi-discrete entropy satisfying for the entropy pair $(\widetilde{\eta}, \widetilde{G})$ in (1.8) is that for some $H^*$ we have*

$$\mathcal{E}_l \geq 0, \qquad \mathcal{E}_r \leq 0, \quad (A.4)$$

*where*

$$\begin{aligned}
\mathcal{E}_l &\equiv \mathcal{F}^h \cdot \left( g\big(h_l - h_l^*\big) - H_l + H^* \right) + \frac{(u_l^*)^2}{2} - \frac{u_l^2}{2} \\
&\quad + (u_l - u_l^*)(\mathcal{F}^q - p(h_l^*)) + u_l T^-, \quad (A.5) \\
\mathcal{E}_r &\equiv \mathcal{F}^h \cdot \left( g\big(h_r - h_r^* - H_r + H^*\big) \right) + \frac{(u_r^*)^2}{2} - \frac{u_r^2}{2} \\
&\quad + (u_r - u_r^*)(\mathcal{F}^q - p(h_r^*)) + u_r T^+. \quad (A.6)
\end{aligned}$$

Consider now the modified hydrostatic reconstruction scheme proposed before and assume we are in Case 2. Assume now that $\mathcal{F}(w_{i+1/2}^-, w_{i+1/2}^+)$ is a conservative flux for the homogeneous problem which satisfies a semi-discrete entropy inequality. Suppose for instance that $h_i - H_i + H_{i+1} < 0$.

If we try to apply Proposition A.2, we get that a sufficient condition for the modified hydrostatic reconstruction scheme to satisfy a semidiscrete entropy inequality is

$$0 \leq \mathcal{F}^h \cdot g(h_i - h_{i+1/2}^- - H_i + H_{i+1}) + u_i T_{i+1/2}^-$$
$$= \mathcal{F}^h \cdot g(h_i - H_i + H_{i+1}) - gh_i u_i(h_i - H_i + H_{i+1}). \quad \text{(A.7)}$$

If $\mathcal{F}$ defines a positive scheme for shallow water with flat bottom, then $\mathcal{F}^h \equiv \mathcal{F}^h(0, w_{i+1/2}^+) \leq 0$ and $\mathcal{F}^h \cdot g(h_i - H_i + H_{i+1}) \geq 0$. Unfortunately, the sign of $-gh_i u_i(h_i - H_i + H_{i+1})$ will depend on the sign of $u_i$.

As a consequence, we may say that the modified hydrostatic reconstruction is semi-discrete entropy satisfying in most cases, but eventually this property could be lost.

## Acknowledgments

This research has been partially supported by the Spanish Government and FEDER through the Research projects MTM2009-11923, and by the Andalusian Government through the project P11- RNM7069. The numerical computations have been performed at the Laboratory of Numerical Methods of the University of Málaga.